\documentstyle[12pt]{article}

\newcommand{\be}{\begin{equation}}
\newcommand{\ee}{\end{equation}}
\newcommand{\ba}{\begin{eqnarray}}
\newcommand{\ea}{\end{eqnarray}}
\newcommand{\ban}{\begin{eqnarray*}}
\newcommand{\ean}{\end{eqnarray*}}
\newcommand{\n}{\nonumber \\}

\newcommand{\eq}[1]{(\ref{#1})}
\newcommand{\sfrac}[2]{{\textstyle \frac{#1}{#2}}}
\newcommand{\ignore}[1]{}

\begin{document}

\renewcommand{\thefootnote}{\fnsymbol{footnote}}
\font\csc=cmcsc10 scaled\magstep1

{\baselineskip=14pt
 \rightline{
 \vbox{
       \hbox{DPSU-96-7}
       \hbox{UT-750}
       \hbox{April 1996}
}}}

\vskip 11mm
\begin{center}

{\large\bf Vertex Operators of the $q$-Virasoro Algebra; \\
\vskip.1in
Defining Relations, Adjoint Actions and Four Point Functions}

\vspace{15mm}

{\csc Hidetoshi AWATA}\footnote{JSPS fellow}$^{1}$,
{\csc Harunobu KUBO}$^{*2}$,
{\csc Yoshifumi MORITA}$^{*3}$,
\\ 
\vskip.05in
{\csc Satoru ODAKE}$^4$ and
{\csc Jun'ichi SHIRAISHI}$^5$.
{\baselineskip=15pt
\it\vskip.35in 
\setcounter{footnote}{0}\renewcommand{\thefootnote}{\arabic{footnote}}
\footnote{e-mail address : hawata@rainbow.uchicago.edu}
James Frank Institute and Enrico Fermi Institute,
University of Chicago,\\
5640, S. Ellis Ave., Chicago, IL, 60637, U.S.A.
\vskip.1in 
\footnote{e-mail address : kubo@danjuro.phys.s.u-tokyo.ac.jp}
Department of Physics, Faculty of Science \\
University of Tokyo, Tokyo 113, Japan \\
\vskip.1in 
\footnote{e-mail address : morita@kodama.issp.u-tokyo.ac.jp}
${}^{\!\!\!,\,5}$ Institute for Solid State Physics, \\
University of Tokyo, Tokyo 106, Japan \\
\vskip.1in 
\footnote{e-mail address : odake@azusa.shinshu-u.ac.jp\\
${}^{\;\;\;\;\;\;5}$e-mail address : shiraish@momo.issp.u-tokyo.ac.jp
}
Department of Physics, Faculty of Science \\
Shinshu University, Matsumoto 390, Japan\\
\vskip.1in 
%
}
\end{center}

\vspace{7mm}
\begin{abstract}
Primary fields of the $q$-deformed Virasoro algebra are constructed.
Commutation relations among the primary fields 
are studied.
Adjoint actions of the deformed Virasoro current on the 
primary fields are represented by the 
shift operator $\Theta_{\xi} f(x)=f(\xi x)$.
Four point functions of the primary fields 
enjoy the connection formula associated with the Boltzmann 
weights of the fusion Andrews-Baxter-Forrester model.
\end{abstract}

\vspace{2mm}
q-alg/9604023

\vfill\eject
\setcounter{footnote}{0}\renewcommand{\thefootnote}{\arabic{footnote}}

{\bf 1. Introduction}\qquad
One of the interesting topics in two dimensional systems
including the conformal field theory (CFT) \cite{rBPZ} and 
solvable lattice models, 
is to develop methods for calculating correlation functions of 
physical observables.
In these solvable models, infinite dimensional symmetries 
play an essential role. 
The most fundamental 
symmetries of the CFT 
are the Virasoro algebra, the Kac-Moody
algebra or their supersymmetric counterparts. 
For the six vertex model and $XXZ$ quantum spin chain, it is 
well recognized that the quantum affine algebra 
$U_q(\widehat{sl}_2)$ takes place.
However, applications of this approach to lattice models had been 
restricted to a class of models which are 
defined by trigonometric solutions of the Yang-Baxter equation.
For some models,
this is because we lack the techniques of bosonizing the 
vertex operators (intertwining operators of the 
symmetry algebras) and for others, 
the characterization of their infinite symmetries are 
still missing.

Recently the symmetry of
the Andrews-Baxter-Forrester (ABF) model \cite{rABF} 
was proposed in the work \cite{rPL,rLP}.
On the other hand, some of the authors found a Virasoro-like
symmetry in a $q$-deformation of the 
Calogero-Sutherland model associated with the 
Macdonald symmetric polynomials \cite{rSKAO} (see also \cite{rFF,rAKOS}). 
A bosonization formula of the deformed Virasoro generator
and the screening operators was also constructed
to study the structure of the highest weight modules of the 
deformed Virasoro algebra. 
It is astonishing that 
this deformed Virasoro algebra was identified with the 
symmetry of the ABF model in \cite{rL,rLP}. 
The important objects in the calculation of the correlation functions
of the ABF model are the vertex operators which can be 
regarded as $q$-deformations of the $(2,1)$ and $(1,2)$ operators 
$\phi_{2,1}(z),\phi_{1,2}(z)$ of
the Virasoro minimal model.

In this letter, we will study a natural way to deform 
$(\ell+1,k+1)$ operators. For this end,
$(\ell+1,1)$ and $(1,\ell+1)$ ($\ell=1,2,\cdots$) operators will be
considered in detail. Using that, some properties of the 
$(\ell+1,k+1)$ operators will be  derived.
One of our requirements for deformed vertex operators 
is that some simple cases of 
their four point functions can be expressed 
by the {\it $q$-hypergeometric functions}, that makes it possible to
have their connection matrices
give the Boltzmann weights of the fusion ABF model.
This is a straightforward extension of the idea in \cite{rLP}.
This assumption, however, does not fix answers to this problem uniquely.
Therefore, another principle must be needed to fix our goal. 
In the CFT, we can define an adjoint action of the 
energy-momentum tensor on the primary fields, and 
that gives us the $c=0$ action of the Virasoro algebra: 
$L_n=-z^{n+1}\frac{d}{dz}$.
We expect that {\it 
we are also able to define an adjoint action of the deformed
Virasoro algebra on the deformed vertex operators}.
So as to obtain this property, it is desirable to have the 
commutation relation between the vertex operator 
$V(z)$ and  
the deformed Virasoro current $T(z)$ as follows:
$g(w/z)T(z)V(w) - V(w)T(z)g(z/w)=\sum_i c_i \delta(\alpha_i w/z) V(\beta_i w)$,
where $g(x)$ is a structure function and $c_i$, $\alpha_i$ and $\beta_i$ 
represent the coefficients, 
the place of the singularity and the shift of the vertex operator respectively.
Our solution for $(\ell+1,1)$ and $(1,\ell+1)$ operator have this
property, but for general $(\ell+1,k+1)$ operators ($\ell,k\geq 1$),
this does not hold.
For these operators, on the other hand, we have commutation relations: 
$g(w/z)T(z)V(w) - V(w)T(z)g(z/w)=
\sum_i c_i \delta(\alpha_i w/z) ( V(w) T(z) \tilde{g}_i(w/z) )$,
where $\tilde{g}(x)$'s are some functions.  

The plan of this letter is the following.
In Section 2, a brief summary of the deformed Virasoro algebra and the 
definition of the vertex operators are given. Commutation relations 
between the deformed Virasoro current and these vertex operators are
studied. A shift operator representation of the deformed Virasoro 
current is derived by studying a deformed adjoint action 
acting on the vertex operators of type 
$(\ell+1,1)$ and $(1,\ell+1)$.
In Section 3, 
four-point functions of the vertex operators with 
one screening operator are calculated explicitly.
In section 4 is devoted to discussion.

After finishing this work, a paper by Kadeishvili \cite{rK} 
appeared. His vertex operators $V_{\ell,k}(z)$ 
are different from ours. While we were seeking 
`good' definitions of the vertex operators,
we also found a similar object as his and 
some algebraic structure of that. In this work, however, 
we will discuss another possibility, 
because we are interested in the
adjoint action of the deformed Virasoro current.

\vspace{1cm}

{\bf 2. Deformed vertex operators}\\
{\sl 2.1. Definition}\quad
First we recall the defining relation of the $q$-Virasoro
algebra \cite{rSKAO}
having two parameters $q$ and $t$. Set 
$\beta=\log t/\log q$ and $p=qt^{-1}$.
The relation is
\ba
  f(\sfrac{w}{z})T(z)T(w)-T(w)T(z)f(\sfrac{z}{w})
  &\!\!=\!\!&
  -(q^{\frac12}-q^{-\frac12})(t^{\frac12}-t^{-\frac12})
  \frac{\delta(p\sfrac{w}{z})-\delta(p^{-1}\sfrac{w}{z})}
  {p^{\frac12}-p^{-\frac12}}, \label{defeq}
\ea
where the structure function $f(x)$ is 
\ba
  f(x) &\!\!=\!\!&
  \exp\biggl(-\sum_{n>0}\frac{1}{n}\frac{(q^{\frac{n}{2}}-q^{-\frac{n}{2}})
  (t^{\frac{n}{2}}-t^{-\frac{n}{2}})}{p^{\frac{n}{2}}+p^{-\frac{n}{2}}}
  x^n\biggr), \nonumber
\ea
and the delta function is defined by $\delta(x)=\sum_{n\in{\bf Z}}x^n$.
The relation (\ref{defeq}) is  invariant under the following transformations,
\be
  {\rm (I)} \qquad 
  (q,t)\rightarrow(q^{-1},t^{-1}),
  \qquad \qquad
   {\rm (II)} \qquad q \leftrightarrow t\quad 
. \label{sym}
\ee
In the following, we respect these symmetries in bosonization formulas.

{\sl 2.2. Bosonization}\quad
Let us introduce the fundamental Heisenberg algebra\footnote{
The bosons $a_n,Q$ in \cite{rSKAO} are related to $h_n,Q_h$ as 
$a_n=-n\frac{1}{1-t^n}h_n$ ($n>0$), 
$a_{-n}=n\frac{1+p^n}{1-t^{-n}}h_{-n}$ ($n>0$),
$a_0=\frac{1}{\sqrt{\beta}}h_0$, $Q=\frac{2}{\sqrt{\beta}}Q_h$, 
and $h^1_n$ in \cite{rAKOS} is $h_n=h^1_np^{-\frac{n}{2}}$.}
$h_n$ ($n\in{\bf Z}$), $Q_{h}$ having the commutation relations
\ba
  [h_n,h_m]
  &\!\!=\!\!&
  \frac{1}{n}\frac{(q^{\frac{n}{2}}-q^{-\frac{n}{2}})
  (t^{\frac{n}{2}}-t^{-\frac{n}{2}})}{p^{\frac{n}{2}}+p^{-\frac{n}{2}}}
   \delta_{n+m,0}
   ,\qquad
  [h_n,Q_h]=\frac12\delta_{n,0}. \label{boscom}
\ea
The symmetries (\ref{sym}) are taken into account by the invariance of 
(\ref{boscom})
under the isomorphisms of the Heisenberg algebra:
\ba
  &&{\rm (I')}\;\; \theta:
 (q,t) \mapsto (q^{-1},t^{-1}),\quad 
  h_n \mapsto  -h_n \:(n\neq 0),
\quad
  h_0 \mapsto h_0, \quad Q_h \mapsto Q_h, \n
  && {\rm (II')} \;\;\omega: \quad 
 q \leftrightarrow t \quad,  \quad 
  h_n \mapsto -h_n,
  \quad Q_h \mapsto - Q_h . \label{symm}
\ea
The $q$-Virasoro generator $T(z)=\Lambda^+(z)+\Lambda^-(z)$ 
and screening currents\footnote{
Screening currents in \cite{rSKAO} (say, $S^{\rm old}_{\pm}(z)$) is 
related to $S_{\pm}(z)$ as
$
  S_+(z)=S_+^{\rm old}(q^{-\frac12}z)q^{\sqrt{\beta}h_0},\;
  S_-(z)=S_-^{\rm old}(t^{-\frac12}z)t^{-\frac{1}{\sqrt{\beta}}h_0}.
$
This modification does not affect the important relations between
singular vectors of the $q$-Virasoro algebra and the Macdonald symmetric
polynomials, because both of them give the same OPE factors.
}  $S_{\pm}(z)$
are bosonized as follows:
\ba
  \Lambda^+(z) &\!\!\!=\!\!\!&
  :\exp\biggl( \sum_{n\neq 0}h_np^{\frac{n}{2}}z^{-n}\biggr):
  q^{\sqrt{\beta}h_0}p^{\frac12}, \quad
  \Lambda^-(z) = \theta\cdot \Lambda^+(z)
  =  \omega\cdot \Lambda^+(z),\\
  S_{+}(z) &\!\!\!=\!\!\!& 
  :\exp\biggl(- \sum_{n\neq 0}
  \frac{p^{\frac{n}{2}}+p^{-\frac{n}{2}}}
 {q^{\frac{n}{2}}-q^{-\frac{n}{2}}}
  h_nz^{-n}\biggr):
  e^{2\sqrt{\beta}Q_h}z^{2\sqrt{\beta}h_0},\quad
  S_{-}(z) = \omega\cdot S_{+}(z),
\ea
where $\omega\cdot \beta$ should be understood as $1/\beta$.
We also have $ \theta\cdot S_{\pm}(z)= S_{\pm}(z)$.
The screening current $S_{+}(z)$  enjoys
the commutation relation
\ba
  \lbrack T(z),S_{+}(w)\rbrack
  &\!\!=\!\!&
  -(1-q)(1-t^{-1})
\frac{d_q}{d_{q}w}
 \biggl(w\:\delta(q^{-\frac12}\sfrac{w}{z})
  :\Lambda^-(q^{-\frac12}w)\:S_{+}(w):\biggr),  \label{scrcom}
\ea
where the difference operator is defined by
$\frac{d_\xi}{d_\xi w}F(w)=\frac{F(w)-F(\xi w)}{(1-\xi)w}$.
A similar difference formula for $S_{-}(z)$ can be derived from 
(\ref{scrcom}) by applying $\omega$. 

{\sl 2.3. Simple vertex operators}\quad
Now we define vertex operators (primary fields) $V_{\ell+1,1}(z)$, 
$V_{1,\ell+1}(z)$ by
\ba
  V_{\ell+1,1}(z) &\!\!=\!\!&
  :\exp\biggl(\sum_{n\neq 0}
  \frac{1}{q^{\frac{n}{2\ell}}-q^{-\frac{n}{2\ell}}}h_nz^{-n}\biggr):
  e^{-\ell\sqrt{\beta}Q_h}z^{-\ell\sqrt{\beta}h_0},
\ea
and $V_{1,\ell+1}(z) = \omega\cdot 
V_{\ell+1,1}(z)$ for $\ell = 1,2,\cdots$. 
They have the invariance $V_{1,\ell+1}(z) = \theta\cdot 
V_{1,\ell+1}(z)$ and $V_{\ell+1,1}(z) = \theta\cdot 
V_{\ell+1,1}(z)$.

The vertex operators
$V_{\ell+1,1}(z)$, $V_{1,\ell+1}(z)$ are expressed as fusion of 
the fundamental ones 
$V_{2,1}(z)$, $V_{1,2}(z)$:
\be
  V_{\ell+1,1}(z)=
  \,:\!\prod_{j=1}^{\ell}V_{2,1}(q^{\frac{\ell+1-2j}{2\ell}}z):,\quad
  V_{1,\ell+1}(z)=
  \,:\!\prod_{j=1}^{\ell}V_{1,2}(t^{\frac{\ell+1-2j}{2\ell}}z):.
\label{fusion}
\ee

{\bf Lemma 2.1.}
{\it We obtain the fundamental commutation relation:\footnote{
We use the standard convention $\prod_{j=n}^{n-1}\ast=1$.}
\ba
  &&
  g^{(\ell+1,1)}(\sfrac{w}{z})\Lambda^{+}(z)V_{\ell+1,1}(w)
  -V_{\ell+1,1}(w)\Lambda^{+}(z)(-\sfrac{z}{w})^{2-\ell}
  g^{(\ell+1,1)}(\sfrac{z}{w}) \n
  &=\!\!& 
  p^{\frac12}t^{-\frac{\ell}{2}}
  \prod_{j=0}^{\ell-2}\Bigl(1-t q^{-\frac{j}{\ell}}\Bigr)\cdot
  \delta(p^{ \frac12}q^{ -\frac{1}{2\ell}}\sfrac{w}{z})
  V_{\ell+1,1}(q^{-\frac{1}{\ell}}w), \label{vocom}
\ea
where 
the structure function $g^{(\ell+1,1)}(x)$ is given by
\ba
  g^{(\ell+1,1)}(x)
  &\!\!=\!\!& 
  \exp\biggl(\sum_{n>0}\frac{1}{n}g^{(\ell+1,1)}_nx^n\biggr), \\
  g^{(\ell+1,1)}_n
  &\!\!=\!\!& 
  -\frac{t^n-t^{-n}}{q^{\frac{n}{2\ell}}-q^{-\frac{n}{2\ell}}}
  \frac{1}{p^{\frac{n}{2}}+p^{-\frac{n}{2}}}
  +\frac{p^{-\frac{n}{2}}q^{\frac{n}{\ell}}-p^{\frac{n}{2}}q^{-\frac{n}{\ell}}}
  {q^{\frac{n}{2\ell}}-q^{-\frac{n}{2\ell}}}.
\ea
}
Note that the structure function $g^{(\ell+1,1)}(x)$ is invariant 
under the transformation (\ref{sym})-$\rm{(I)}$.
To prove (\ref{vocom}), the identity
$$
  e^{\sum_{n>0}\frac{1}{n}(1-r_1^n-\cdots-r_m^n)x^n}
  -(-x)^{m-1}\prod_{i=1}^mr_i\cdot 
  e^{\sum_{n>0}\frac{1}{n}(1-r_1^{-n}-\cdots-r_m^{-n})x^{-n}}
  =
  \prod_{i=1}^m(1-r_i)\cdot\delta(x), 
$$
may be useful.
The commutation relations among $T(z)$ and $ V_{\ell+1,1}(z)$,
$ V_{1,\ell+1}(z)$ are derived from (\ref{vocom}) by applying 
$\theta$ and $\omega$.

{\bf Proposition 2.1.}
{\it For $ V_{\ell+1,1}(z)$, we have
\ba
  &&
  g^{(\ell+1,1)}(\sfrac{w}{z})T(z)V_{\ell+1,1}(w)
  -V_{\ell+1,1}(w)T(z)(-\sfrac{z}{w})^{2-\ell}
  g^{(\ell+1,1)}(\sfrac{z}{w}) \n
  &=\!\!& 
  \prod_{j=0}^{\ell-2}
  \Bigl( t^{- \frac12} q^{\frac{j}{2\ell}}-
     t^{ \frac12} q^{-\frac{j}{2\ell}}\Bigr) 
  \left(p^{ \frac12} t^{ -\frac12} q^{ -\frac{(\ell-1)(\ell-2)}{4\ell}} 
  \delta(p^{ \frac12}q^{ -\frac{1}{2\ell}}\sfrac{w}{z})
 V_{\ell+1,1}(q^{-\frac{1}{\ell}}w) \right. \\
&&\qquad\qquad\qquad\quad
-\left. (-1)^{\ell} p^{ -\frac12} t^{ \frac12}
  q^{\frac{(\ell-1)(\ell-2)}{4\ell}} 
  \delta(p^{ -\frac12}q^{ \frac{1}{2\ell}}\sfrac{w}{z})
 V_{\ell+1,1}(q^{\frac{1}{\ell}}w)
\right). \nonumber
\ea
}

Let us define the adjoint action of the $q$-Virasoro generator 
on the vertex operator by
\begin{eqnarray}
 {\cal T}_n^{\ell} \cdot V_{\ell+1,1}(w)
&=&
\oint \frac{dz}{2 \pi iz} z^n \left(
 g^{(\ell+1,1)}(\sfrac{w}{z})T(z)V_{\ell+1,1}(w) \right. \n
  &&\qquad\qquad
-\left. V_{\ell+1,1}(w)T(z)(-\sfrac{z}{w})^{2-\ell}
  g^{(\ell+1,1)}(\sfrac{z}{w}) \right).
\end{eqnarray}
Then we obtain

{\bf Theorem 2.1.}
{\it The operator ${\cal T}_n^{\ell}$ can be represented by the
shift operator $ \Theta_{\xi}$ defined by $ \Theta_{\xi} f(z)=f(\xi
z)$ as
\begin{eqnarray}
{\cal T}_n^{\ell} &=&
 \prod_{j=0}^{\ell-2}
  \Bigl( t^{- \frac12} q^{\frac{j}{2\ell}}-
     t^{ \frac12} q^{-\frac{j}{2\ell}}\Bigr)\cdot
\left(
p^{ \frac{n+1}{2}} t^{-\frac12}
q^{ -\frac{(\ell-1)(\ell-2)+2n}{4\ell}}w^n 
 \Theta_{q^{-\frac{1}{\ell}}} \right. \n
&&\qquad\qquad\qquad\qquad
-(-1)^{\ell}\left.
p^{-\frac{n+1}{2} } t^{\frac12}
q^{ \frac{(\ell-1)(\ell-2)+2n}{4\ell}   } w^n
 \Theta_{ q^{\frac{1}{\ell}}} \right),
\end{eqnarray}
on the vertex operator $V_{\ell+1,1}(w)$.}

{\sl 2.4. General vertex operators}\quad
The following lemma is helpful to construct general 
vertex operators of type $(\ell+1,k+1)$.

{\bf Lemma 2.2.}
{\it 
 The fundamental relation (\ref{vocom}) can be written  in another way,
\ba
  &&
  g^{(\ell+1,1)}(\sfrac{w}{z})\Lambda^{+}(z)V_{\ell+1,1}(w)
  -V_{\ell+1,1}(w)\Lambda^{+}(z)(-\sfrac{z}{w})^{2-\ell}
  g^{(\ell+1,1)}(\sfrac{z}{w})  \n
  &=\!\!& 
  t^{-\frac{\ell}{2}}
  \prod_{j=0}^{\ell-2}\Bigl(1-t q^{-\frac{j}{\ell}}\Bigr)\cdot
  \delta(p^{ \frac12}q^{ -\frac{1}{2\ell}} \sfrac{z}{w})  
\left(
\tilde{g}^{(\ell+1,1)}(\sfrac{w}{z}) V_{\ell+1,1}(w) T(z)\right),\label{VTcom}
\ea
where 
\begin{equation}
\tilde{g}^{(\ell+1,1)}(x) 
=
  \exp\left( - \sum_{n>0}
  \frac{1}{n}\frac{ p^{-\frac{n}{2}}
  (q^{\frac{n}{2}}-q^{-\frac{n}{2}})
  (t^{\frac{n}{2}}-t^{-\frac{n}{2}}) }
  { (q^{\frac{n}{2\ell}} -q^{-\frac{n}{2\ell}} )
  (p^{\frac{n}{2}} +p^{-\frac{n}{2} })}
  x^n   
  \right).
\ee}
To obtain (\ref{VTcom}) from (\ref{vocom}), we used the identity
$V_{\ell+1,1}(q^{-\frac{1}{\ell}}w)=p^{-\frac12}
:V_{\ell+1,1}(w) \Lambda^{+}(p^{\frac12}q^{-\frac{1}{2\ell}}w):$.
It may be worth to make a comment on the 
R.H.S. The factor $V_{\ell+1,1}(w)\Lambda^{-}(z)$ is regular at 
$z=p^{ \frac12}q^{ -\frac{1}{2\ell}} w$, and 
$V_{\ell+1,1}(w)\Lambda^{+}(z)$
has a simple pole at the point. However, the function
$\tilde{g}^{(\ell+1,1)}(w/z) $
has a simple zero at the point. So, the R.H.S. is well defined in totality.
\ignore{
\ba
  &&
  g^{(\ell+1,1)}(\sfrac{w}{z})\Lambda^{+}(z)V_{\ell+1,1}(w)
  -V_{\ell+1,1}(w)\Lambda^{+}(z)(-\sfrac{z}{w})^{2-\ell}
  g^{(\ell+1,1)}(\sfrac{z}{w})  \n
  &=\!\!& 
  p^{\frac12}t^{-\frac{\ell}{2}}
  \prod_{j=0}^{\ell-1}\Bigl(1-t q^{-\frac{j}{\ell}}\Bigr)\cdot
  \delta(p^{ \frac12}q^{ -\frac{1}{2\ell}} \sfrac{w}{z})   \n
 &&  
\left\{ 
\exp\left\{  - \sum_{n>0}
\frac{1}{n}\frac{ p^{-\frac{n}{2}}
 (q^{\frac{n}{2}}-q^{-\frac{n}{2}})
 (t^{\frac{n}{2}}-t^{-\frac{n}{2}}) }
{ (q^{\frac{n}{2\ell}} -q^{-\frac{n}{2\ell}} )
  (p^{\frac{n}{2}} +p^{-\frac{n}{2} })}
\left( \frac{w}{z} \right)^n   
\right\}
V_{\ell+1,1}(w) T(z)p^{-1/2}
\right\}
\ea

\be
V_{\ell+1,1}(q^{-\frac{1}{\ell}}w)
=:V_{\ell+1,1}(w) \Lambda^{+}(p^{\frac12} q^{-\frac{1}{2\ell}}w): p^{-\frac12}
\ee
\ba
&&
  \delta(p^{ \frac12}q^{ -\frac{1}{2\ell}} \sfrac{w}{z})
   V_{\ell+1,1}(q^{-\frac{1}{\ell}}w)  \n
&=& 
\delta(p^{ \frac12}q^{ -\frac{1}{2\ell}} \sfrac{w}{z})
V_{\ell+1,1}(w)\Lambda^{+}(z) p^{-\frac12} 
\exp
\left\{  - \sum_{n>0}
\frac{1}{n}\frac{ p^{-\frac{n}{2}}
 (q^{\frac{n}{2}}-q^{-\frac{n}{2}})
 (t^{\frac{n}{2}}-t^{-\frac{n}{2}}) }
{ (q^{\frac{n}{2\ell}} -q^{-\frac{n}{2\ell}} )
  (p^{\frac{n}{2}} +p^{-\frac{n}{2} })}
\left( \frac{w}{z} \right)^n   
\right\}
\ea
}

Let us define the vertex operator of type $(\ell+1,k+1)$ by
\be
  V_{\ell+1,k+1}(z)=\,:V_{\ell+1,1}(z)V_{1,k+1}(z):,
\ee
and introduce the structure functions: 
$g^{(1,k+1)}(x)= \omega\cdot g^{(k+1,1)}(x)$, 
$
  \tilde{g}^{(1,k+1)}(x)= \theta\cdot\omega\cdot\tilde{g}^{(k+1,1)}(x),
$
$
  g^{(\ell+1,k+1)}(x)=g^{(\ell+1,1)}(x)g^{(1,k+1)}(x),
$
and 
$
  \tilde{g}^{(\ell+1,k+1)}(x)=
  \tilde{g}^{(\ell+1,1)}(x)\tilde{g}^{(1,k+1)}(x).
$
Using Lemma 2.2. and the maps $\theta$ and $\omega$, we have the
commutation relation for this vertex operator.

{\bf Proposition 2.2.}
{\it The commutation relation between the deformed Virasoro current
and the vertex operator of type $(\ell+1,k+1)$ 
is given by the following relation and $\theta$, $\omega$ operations,}
\ba
  &&
  g^{(\ell+1,k+1)}(\sfrac{w}{z})\Lambda^{+}(z)V_{\ell+1,k+1}(w)
  -V_{\ell+1,k+1}(w)\Lambda^{+}(z)(-\sfrac{z}{w})^{4-\ell-k}
  g^{(\ell+1,k+1)}(\sfrac{z}{w}) \n
  &=\!\!& 
  V_{\ell+1,k+1}(w)T(z)\tilde{g}^{(\ell+1,k+1)}(\sfrac{z}{w})
  q^{\frac{k}{2}}t^{-\frac{\ell}{2}} \n
  &&
  \times\biggl(
  \frac{\delta(p^{\frac12}q^{-\frac{1}{2\ell}}\sfrac{w}{z})}
  {1-q^{\frac{1}{2\ell}}t^{\frac{1}{2k}}}
  \prod_{i=0}^{\ell-2}(1-tq^{-\frac{i}{\ell}})\cdot
  \prod_{j=0}^{k-2}(1-q^{\frac{1}{2\ell}-1}t^{\frac{1}{2k}+\frac{j}{k}})
  \label{TVkl} \\
  &&\qquad
  +\frac{\delta(p^{\frac12}t^{\frac{1}{2k}}\sfrac{w}{z})}
   {1-q^{-\frac{1}{2\ell}}t^{-\frac{1}{2k}}}
  \prod_{j=0}^{k-2}(1-q^{-1}t^{\frac{j}{k}})\cdot
  \prod_{i=0}^{\ell-2}
   (1-t^{1-\frac{1}{2k}}q^{-\frac{1}{2\ell}-\frac{i}{\ell}})\biggr).
  \nonumber
\ea

In the $q\rightarrow 1$ limit, $V_{r,s}(z)$ reduces to the usual vertex
operator $e^{\frac12\alpha_{r,s}\phi(z)}$  
where $\alpha_{r,s}=\sqrt{\beta}(1-r)-\frac{1}{\sqrt{\beta}}(1-s)$ 
and $\phi(z)\phi(w)\sim 2\log(z-w)$. 
This can be easily shown in the bosonized form.
(We remark that there are infinitely many operators which reduce to
$e^{\frac12\alpha_{r,s}\phi(z)}$ in the $q\rightarrow 1$ limit.)  
However, it is rather nontrivial
to derive the usual defining relation of the
primaly field of the Virasoro algebra from the commutation relations 
\eq{TVkl}.

\vspace{1cm}

{\bf 3. Correlation functions}\qquad 
First we recall standard notations;
$(a;q_1,\cdots,q_{\ell})_n
 =\prod_{k_1,\cdots,k_{\ell}=0}^{n-1}
 (1-aq_1^{k_1}\cdots q_{\ell}^{k_{\ell}})$, 
$(a_1,\cdots,a_m;q_1,\cdots,q_{\ell})_n
  =\prod_{j=1}^m(a_j;q_1,\cdots,q_{\ell})_n$,  
$\Gamma_q(z)=(1-q)^{1-z}(q;q)_{\infty}/(q^z;q)_{\infty}$,
$B_q(x,y)=\Gamma_q(x)\Gamma_q(y)/\Gamma_q(x+y)$,
$\vartheta_q(z)=(q;q)_{\infty}(z;q)_{\infty}(qz^{-1};q)_{\infty}$.
The Jackson integral is defined by 
$\int_0^ad_qzf(z)\!\!=\!\!a(1-q)\sum_{n=0}^{\infty}f(aq^n)q^n$,  
$\int_0^{a\infty}d_qzf(z)\!\!=\!\!a(1-q)\sum_{n=-\infty}^{\infty}f(aq^n)q^n$, 
$\int_A^Bd_qzf(z)\!\!=\!\!\int_0^Bd_qzf(z)-\int_0^Ad_qzf(z)$.
The $q$-hypergeometric function ${}_2\phi_1$ is defined by
\be
  {}_2\phi_1(a,b;c;q,z)
  =\sum_{n=0}^{\infty}\frac{(a;q)_n(b;q)_n}{(q;q)_n(c;q)_n}z^n.
\ee

In the following we abbreviate $V_{\ell+1,1}(z)$ as $V_{\ell}(z)$.
We will calculate the following four point functions,
\ba
  U_+(z,w)&\!\!=\!\!&
  \int_{t^{\frac12}q^{\frac{1}{2\ell}}z}^{
        t^{\frac12}q^{\frac{1}{2\ell}}z\infty}
  d_{q^{\frac{1}{\ell}}}\mu
  \langle\ast|S_+(\mu)V_{\ell}(z)V_{\ell}(w)V_L(0)|0\rangle, \\
  U_-(z,w)&\!\!=\!\!&
  \int_0^{t^{-\frac12}q^{\frac{1}{2\ell}}w}
  d_{q^{\frac{1}{\ell}}}\mu
  \langle\ast|V_{\ell}(z)V_{\ell}(w)S_+(\mu)V_L(0)|0\rangle,
\ea
where the momentum of the bra state $\ast$ is chosen such that 
$\langle\ast|:S_+(\mu)V_{\ell}(z)V_{\ell}(w)V_L(0):|0\rangle=1$.
These four point functions are expressed as
\ba
  && U_+(z,w)/\langle V_{\ell}(z)V_{\ell}(w)\rangle \n
  &=\!\!&
  (zw)^{\frac12\ell L\beta} 
  B_{q^{\frac{1}{\ell}}}(2a-b,1-a)(t^{\frac12}q^{-\frac{1}{2\ell}}z)^{b-2a}
  {}_2\phi_1(q^{\frac{a}{\ell}},q^{\frac{2a-b}{\ell}};q^{\frac{1+a-b}{\ell}};
  q^{\frac{1}{\ell}},q^{\frac{1-a}{\ell}}\sfrac{w}{z}), \\
  &&  U_-(z,w)/\langle V_{\ell}(z)V_{\ell}(w)\rangle \n
  &=\!\!& 
  (zw)^{\frac12\ell L\beta} 
  B_{q^{\frac{1}{\ell}}}(b,1-a)(zw)^{-a}(t^{-\frac12}q^{\frac{1}{2\ell}}w)^b
  {}_2\phi_1(q^{\frac{a}{\ell}},q^{\frac{b}{\ell}};q^{\frac{1-a+b}{\ell}};
  q^{\frac{1}{\ell}},q^{\frac{1-a}{\ell}}\sfrac{w}{z}),
\ea
where $ a=\ell\beta$, $b=1-L\beta$ and 
\begin{equation}
  \langle V_{\ell}(z)V_{\ell}(w)\rangle
  =
  z^{\ell^2 \beta /2}\prod _{j=1}^{\ell}
  \frac{(t^{-1}q^{j/\ell}w/z,pq^{j/\ell}w/z;p^2,q^{1/\ell})_{\infty}}
     {(q^{j/\ell}w/z,t^{-1}pq^{j/\ell}w/z;p^2,q^{1/\ell})_{\infty}}.
\end{equation}
We have used the formulas,
\begin{eqnarray}
S_{+}(\mu)V_{\ell}(z)&\!\!=\!\!&
:S_{+}(\mu)V_{\ell}(z):
\mu^{-\ell \beta}\frac{(t^{1/2} q^{1/2\ell}z/\mu;q^{1/\ell})_{\infty}}
{(t^{-1/2} q^{1/2\ell}z/\mu;q^{1/\ell})_{\infty}},\label{sv}\\
V_{\ell}(z)S_{+}(\mu)
&\!\!=\!\!&
:V_{\ell}(z)S_{+}(\mu):
z^{-\ell\beta}
\frac
{(t^{1/2}q^{1/2\ell}\mu/z;q^{1/\ell})_{\infty}}
{(t^{-1/2}q^{1/2\ell}\mu /z;q^{1/\ell})_{\infty }}.\label{vs}
\end{eqnarray}
We note that
the ratio of the coefficients in the R.H.S.'s of (\ref{sv}) and
(\ref{vs}) is a pseudo-constant
with respect to the shift $z/\mu\rightarrow q^{1/\ell}z/\mu$.

Introduce the notation $x=q^{\frac{1}{2r\ell}}$ and 
the definition 
$[u]\!=\!\sqrt{2 \pi r/\varepsilon} x^{r/4}
x^{u(u-r)/r}\vartheta_{x^{2r}}(x^{2u})$
where $\varepsilon = -2 \pi^2 /\log x$.

{\bf Proposition 3.1.}
{\it
The connection formula for the four point functions $U_{\pm}(z,w)$ 
can be written as 
\ba
  &\!\!\!\!\!\!&
  \left(\begin{array}{c}U_+(z,w)\\U_-(z,w)\end{array}\right)
  =\frac{\langle V_{\ell}(z)V_{\ell}(w)\rangle }
        {\langle V_{\ell}(w)V_{\ell}(z)\rangle }
    \left(\begin{array}{cc}
  \frac{[\ell][-u+\ell+L]}{[\ell+L][u+\ell]} &
  \frac{[L][-u]}{[\ell+L][u+\ell]} \\
  \frac{[2\ell+L][-u]}{[\ell+L][u+\ell]} &
  \frac{[\ell][u+\ell+L]}{[\ell+L][u+\ell]}
  \end{array}\right)  
  \left(\begin{array}{c}U_+(w,z)\\U_-(w,z)\end{array}\right),
  \label{conn-mat}
\ea
where $u$ is defined by $\frac{w}{z}=x^{2u}$.}

The elements of the connection matrix are identified with 
some of the Boltzmann weights of 
the $\ell\times\ell$ fusion RSOS model \cite{rABF}.
We conjecture that 
the connection matrices
for the four point functions having 
arbitrary numbers of 
screening operators are expressible in the same way 
by the Boltzmann weights of the $\ell\times\ell$ fusion ABF model.
\vspace{1cm}

{\bf 4. Discussion}\qquad

We constructed the vertex operators (primary fields) 
of the deformed Virasoro algebra.
Using the vertex operators, we 
obtained the adjoint action of the deformed Virasoro 
current acting on the vertex operator of type $(\ell+1,1)$, 
$(1,\ell+1)$. This kind of adjoint action 
is found only for this simple cases, so far. 
Thus
it is desirable
having more general definition of this kind of adjoint action.

Applications of the vertex operators
might be possible to 
integrable systems; integrable massive field theories, 
solvable lattice models and so on.
One of these candidates is 
the calculation of the correlation functions
of the fusion ABF model, since 
we get some of the Boltzmann weights 
of the $\ell\times\ell$ fusion ABF model
in the connection formula of the four point functions
derived in Section 3. 
As is discussed by Lukyanov and Pugai \cite{rPL,
rLP},
the vertex operators
$V_{2,1}(z)$ and $V_{1,2}(z)$
act on the physical space of the ABF model.
Then, it may be natural to consider that our vertex operators 
$V_{\ell +1,1}(z)$ and $V_{1,\ell +1}(z)
(\ell \ge 2)$
act on some fusion ABF model.
The situation is, however, slightly complicated
i.e.
the 
central charge 
of this model in the regrme III,
which is derived from the corner transfer matrix method, 
is greater than one 
except for the special case: $\# {\rm of\;states} = \# {\rm of\; fusion} +2$.
Thus it seems impossible to represent the physical space of the 
general fusion ABF model 
by a Fock space of single bosonic field.
Therefore, the first task should be to 
identify the model which 
our vertex operators are associated with, if there is any.
To this end,
it might be suggestive
to study 
the factor 
$\frac{\langle V_{\ell}(z)V_{\ell}(w)\rangle }
{\langle V_{\ell}(w)V_{\ell}(z)\rangle }$,
which appears in 
the connection matrix (\ref{conn-mat}).
In the case $\ell =1$, it reduces to 
$\frac {[u+1]/[1]}{\kappa (u)},$ 
where
$\kappa (u)$ 
is the free eneregy of the ABF model in the regeme III.

In the definition of the correlation functions given in Section 3,
the cycle of these Jackson integrals
are given by hand 
depending on 
the ordering of the vertex operators and screening operators.
In the conformal field theory,
however,
the screened vertex operator is defined by
the contour integral 
whose contour is independent of this kind of ordering \cite {rF}. 
Thus,
it would be desirable
to obtain the screened vertex operator
by a contour integral \cite{rPL,rLP,rK}.

\vskip5mm
\noindent{\bf Acknowledgments}\\
We would like to thank
T.~Eguchi,  A.~Kuniba, S.~Lukyanov, Y.~Pugai and T.~Takagi
for valuable discussion.



\begin{thebibliography}{99}

\bibitem{rBPZ}
  A.A.~Belavin, A.M.~Polyakov and A.B.~Zamolodchikov,
    {\sl Nucl. Phys.} {\bf B241} (1984) 333, \\
  {\it Conformal Invariance and
       Applications in Statistical Mechanics,}
  edited by C.~Itzykson, H.~Saleur, and J.-B.~Zuber,
  {\sl (World Scientific, Singapore, 1988)}. 

\bibitem{rABF}
  G.E.~Andrews, R.J.~Baxter and P.J.~Forrester,
  {\it Eight-Vertex SOS Model
   and Generalized Rogers-Ramnujan Identities,}
  {\sl J. Stat. Phys.} {\bf 35} (1984) 193-266,\\
  E.~Date, M.~Jimbo, A.~Kuniba, T.~Miwa and M.~Okado,
  {Exactly Solvable SOS Models: Local height probabilities and theta
    function identities},
  {\sl Nucl. Phys.} {\bf B 290} (1987) 231-273;
  {\it Exactly Solvable SOS Models 2:
   Proof of the Star-Triangle Relation and
   Combinatorial Identities.}
  {\sl Adv. Stud. Pure. Math.} {\bf 16} (1988) 17-122,\\
  {\it Yang-Baxter Equation in Integrable Systems},
  edited by M.~Jimbo,  
  {\sl (World Scientific, Singapore, 1989)}

\bibitem{rPL}
  {S.~Lukyanov and  Y.~Pugai, 
  {\it Bosonization of ZF Algebras:
   Direction Toward Deformed Virasoro Algebra},
  hepth/9412128.}

\bibitem{rSKAO}
  J.~Shiraishi, H.~Kubo, H.~Awata and S.~Odake,
  {\it A Quantum Deformation of the Virasoro Algebra and
  the Macdonald Symmetric Functions},
  q-alg/9507034, to appear in {\sl Lett. Math. Phys.}

\bibitem{rAKOS}
  H.~Awata, H.~Kubo, S.~Odake and J.~Shiraishi,
  {\it Quantum ${\cal W}_N$ Algebras and Macdonald Polynomials},
  q-alg/9508011, to appear in {\sl Comm. Math. Phys.}

\bibitem{rFF}
  B.~Feigin and E.~Frenkel,
  {\it Quantum ${\cal W}$-algebras and elliptic algebras},
  q-alg/9508009.

\bibitem{rL}
  S.~Lukyanov,
  {\it A Note on the Deformed Virasoro Algebra},
  {\sl Phys. Lett.} {\bf B367} (1996) 121-125.

\bibitem{rLP}
  S.~Lukyanov and Y.~Pugai,
  {\it Multi-point Local Height Probabilities in the Integrable
  RSOS Model},
  hep-th/9602074.

\bibitem{rK}
  A.~A.~Kadeishvili, 
  {\it Vertex Operators for Deformed Virasoro Algebra},
 hepth/9604153.


\bibitem{rGR}
  G.~Gasper and M.~Rahman,
  {\it Basic Hypergeometric Series},
  {\sl (Cambridge University Press, 1990)}.
\bibitem{rF}
  G.~Felder:
  {\it BRST Approach to Minimal Models,}
  {\sl Nucl. Phys.}
  {\bf B317} (1989) 215. 

\end{thebibliography}
\end{document}